\documentclass[10pt]{article}
 \usepackage{subeqnarray}
\usepackage{color}
\begin{document}

\newcommand{\nc}{\newcommand}

\newtheorem{lemma}{Lemma}[section]
\newtheorem{theorem}[lemma]{Theorem}
\newtheorem{proposition}[lemma]{Proposition}
\newtheorem{corollary}[lemma]{Corollary}
\newtheorem{remark}[lemma]{Remark}
\newtheorem{example}[lemma]{Example}
\newtheorem{hypothesis}[lemma]{Hypothesis}
\newtheorem{notation}[lemma]{Notation}
\newtheorem{definition}[lemma]{Definition}
\newtheorem{conclusion}[lemma]{Conclusion}

\nc{\QED}{\mbox{}\hfill \raisebox{-0.2pt}{\rule{5.6pt}{6pt}\rule{0pt}{0pt}}
          \medskip\par}
\newenvironment{Proof}{\noindent
    \parindent=0pt\abovedisplayskip = 0.5\abovedisplayskip
    \belowdisplayskip=\abovedisplayskip{\bf Proof. }}{\QED}
\newenvironment{Proofof}[1]{\noindent
    \parindent=0pt\abovedisplayskip = 0.5\abovedisplayskip
    \belowdisplayskip=\abovedisplayskip{\bf Proof of #1. }}{\QED}
\newenvironment{Example}{\begin{example}
                \parindent=0pt \rm}{\QED\end{example}}
\newenvironment{Remark}{\begin{remark}
                \parindent=0pt\rm }{\QED\end{remark}}
\newenvironment{Definition}{\begin{definition}
                \parindent=0pt \rm}{\QED\end{definition}}
\newenvironment{Conclusion}{\begin{conclusion}
                \parindent=0pt \rm}{\QED\end{conclusion}}

\def\ba{\begin{array}}
\def\ea{\end{array}}

\def\be{\begin{equation}}
\def\ee{\end{equation}}
\def\vs5{\vspace{0.5cm}}
\def\lab{\label}
\def\bthm{\begin{theorem}}
\def\ethm{\end{theorem}}
\def\bhyp{\begin{hypothesis}}
\def\ehyp{\end{hypothesis}}
\def\bP{\begin{Proof}}
\def\eP{\end{Proof}}
\def\bPof{\begin{Proofof}}
\def\ePof{\end{Proofof}}
\def\brem{\begin{remark}}
\def\erem{\end{remark}}
\def\bex{\begin{example}}
\def\eex{\end{example}}
\def\bcor{\begin{corollary}}
\def\ecor{\end{corollary}}
\def\bd{\begin{definition}}
\def\ed{\end{definition}}
\def\bprop{\begin{proposition}}
\def\eprop{\end{proposition}}
\def\blem{\begin{lemma}}
\def\elem{\end{lemma}}
\def\beq{\begin{eqnarray}}
\def\eeq{\end{eqnarray}}
\def\beqs{\begin{eqnarray*}}
\def\eeqs{\end{eqnarray*}}

\def\a{\alpha}
\def\b{\beta}
\def\s{\sigma}
\def\Sig{\Sigma}
\def\d{\delta}
\def\l{\lambda}
\def\hs{\hat{\sigma}}
\def\hT{\hat{T}}
\def\hU{\hat{U}}
\def\eps{\epsilon}
\def\veps{\varepsilon}
\def\benum{\begin{enumerate}}
\def\eenum{\end{enumerate}}
\def\bit{\begin{itemize}}
\def\eit{\end{itemize}}
\def\la{\langle}
\def\ra{\rangle}
\def\wto{\rightharpoonup}
\def\bwto{\buildrel{w(\mu_h,\mu)}\over\longrightarrow}
\def\bYto{\buildrel{Y(\mu_h,\mu)}\over\longrightarrow}
\def\bsto{\buildrel{s(\mu_h,\mu)}\over\longrightarrow}
\def\bswto{\buildrel{sw(\mu_h,\mu)}\over\longrightarrow}

\def\Cinf{C^\infty}
\def\Linf{L^\infty}
\def\supp{{\rm supp}}
\def\setm{\setminus}
\def\id{{\rm id}}

\def\dis{\displaystyle}

\def\R{{\bf R}}
\def\N{{\rm N}}
\def\Ra{{\rm R}}
\def\M{{\cal M}}
\def\B{{\cal B}}
\def\A{{\cal A}}
\def\C{{\cal C}}
\def\D{{\bf D}}
\def\calS{{\cal S}}
\def\O{{\cal O}}
\def\P{{\cal P}}
\def\Q{{\cal Q}}
\def\u{{\cal U}}
\def\V{{\cal V}}
\def\F{{\cal F}}
\def\G{{\cal G}}
\def\H{{\cal H}}
\def\E{{\cal E}}
\def\m{{\bf m}}
\def\o{{\bf o}}
\def\LL{{\cal L}}
\def\MM{{\bf M}}
\def\DD{{\cal D}}
\def\RR{{\cal R}}
\def\NN{{\bf N}}
\def\EE{{\bf E}}
\def\FF{{\bf F}}
\def\GG{{\bf G}}
\def\T{{\cal T}}
\def\TT{{\bf T}}
\def\BB{{\bf B}}
\def\HH{{\bf H}}
\def\KP{{\bf KP}}
\def\K{{\bf K}}
\def\aa{{\bf a}}

\def\rank{{\rm Rank}}
\def\span{{\rm span}}
\def\dim{{\rm dim}}
\def\diam{{\rm diam}}
\def\diag{{\rm diag}}
\def\dist{{\rm dist}}
\def\trace{{\rm trace}}
\def\div{{\rm div}}
\def\spt{{\rm spt}}
\def\Lip{{\rm Lip}}
\def\Sp{{\bf S}}
\def\Z{{\bf Z}}

\def\Grass{{\rm Grass}}
\def\diam{{\rm diam}\,}
\def\rmB{{\rm B}}

\def\pa{\partial}
\def\ov{\overline}
\def\na{\nabla}
\def\vpi{\varpi}

\def\ti{\tilde}
\def\ess{{\rm ess}}
\def\ext{^{\rm ext}}
\def\u{{\bf u}}

\def\w{{\bf w}}
\def\e{{\bf e}}
\def\Om{\Omega}
\def\om{\omega}
\def\al{\aleph}

\def\ioe{\int_{\Omega_\veps}}
\def\ioesmbe{\int_{\Omega_\veps\setminus \C_\veps}}
\def\io{\int_{\Omega}}
\def\ise{\int_{D_\veps}}
\def\ibe{\int_{\C\veps}}
\def\iyke{\int_{Y^k_\veps}}
\def\dme{\;dm_\veps}
\def\Ome{{\Omega_\veps}}
\def\me{{m_{\veps}}}
\def\intb{{\int\!\!\!\!\!\!-}}

\def\QT{{\Omega^T}}
\def\De{{D_\veps}}
\def\Rex{{R_\veps}}
\def\Ri{{r_\veps}}

\def\card{{\rm card}}

\title{\textbf{Homogenizing  media containing a  highly conductive honeycomb substructure
}}

\author{\textbf{ by $\,$Isabelle~Gruais $^{*}$ and  $\,$Dan~Poli\v{s}evski $^{**}$
} }
\date{}
\maketitle

{\bf Abstract.} The present paper deals with the homogenization of the heat conduction which takes place in a binary three-dimensional medium  consisting of an ambiental  phase having  conductivity of unity order and a rectangular honeycomb structure formed  by a set of thin layers crossing orthogonally and periodically. We consider the case when the conductivity of the thin layers  is in inverse proportion to the vanishing volume of the rectangular honeycomb structure. We find the system that governs  the asymptotic behaviour of the temperature distribution of this binary medium. The dependence with respect to the thicknesses of the layers is also emphasized. We use an energetic method associated to a natural control-zone of the vanishing domain.



\section{Introduction}\lab{s:1}

\hspace{0.5cm}

The study of the lattice-type structures which are characterized by periodicity and small thickness of the material is one of the main achievements of the homogenization theory. Civil engineering, electrotechnics and the aerospace industry are mainly concerned by composite materials, particularly truss structures,  which have to be treated by techniques of homogenization   when  direct computations fail.

The foundations of homogenized  layered materials were laid down by Murat and Tartar in their
pioneering work \cite{Murattartar85}. This method was  still  used  by \cite{HeronMP94} in the framework of weaker topologies. The case of $BV$-functions and sequences of measures is worked out in \cite{BouchitteP96}. Many examples and applications  may be found in \cite{BakhvalovP89}.  The difficulties arising from the direct computation of the behaviour of these structures are twofold: the great number of cells and the small thickness of the material. The periodic distribution classically suggests that the homogenization method for perforated domains should  be used  as in \cite{SJP99},  where the period $\varepsilon$  and the thickness $\delta$ of the material are considered as independent  vanishing parameters. A classical  issue is then to consider the problem of permuting both converging processes. This is made in details in \cite{SJP99}. As the pratical point of view favors geometric consideration, measures provide a more general tool when they are confined to the description of the critical part of a system,  as in \cite{BellieudB} and \cite{BGP06}. This approach was used in \cite{ZhikovP06}, where the homogenization of elastic reticulated structures is performed for an anisotropic material surrounded by an empty environment.

Our study is based on the control-zone method which was introduced and developed by \cite{BGP}-\cite{BGPprep} and \cite{GP07}-\cite{GP08},  specific  to the binary  composition of the system, namely  the truss and the ambiental  phase. This procedure   proved its efficiency in  the modelling of fine substructures where concentrated material  of high contribution  influences  the behavior of the global problem in spite of  its  vanishing volume. The asymptotic treatment  reveals the apparent paradox
between an  obviously disappearing element  and its everlasting
action on their environment. The coupling between both components manifests itself through the rarefying  ratio $\gamma_\veps$ which arises as a criterium   for the reduced problem to exist. Unlike the critical case studied in  \cite{BGPprep}, the connections  between the layers of the reticulated structure  annihilate the capacitary  term $\lim_{\veps\to 0}\gamma_\veps$  when it is defined. This was also observed in \cite{GP08} in the absence of connections. Interestingly, the presence of the ambiental phase substantially modifies the contribution of the truss in a way that cannot be anticipated from \cite{SJP99}-\cite{ZhikovP06}, because  the sequence of problems under consideration behaves singularly with respect to the period of distribution and involves a new criterium, namely the capacity of the intersections of the crossing layers.

The present paper is organized as follows:

Section~\ref{s:2} is devoted to
the main notations and to the description of the initial problem. We
set the functional framework  for which  the existence and
uniqueness of the solution can be established.

 In the first part  of Section~\ref{s:3}  we  present  two  operators that have a localizing effect  and that we use in order to obtain the specific  inequalities related to honeycomb structures. We also define   the capacitary functions and step approximation operators which are associated with  the test-functions  defined in our control-zones. They overcome  the singular behavior of the energy term when  $\veps$ tends to zero.

Section~\ref{s:4}   deals  with the homogenization process in the  reticulated  case in the so-called box structure geometry. We obtain the homogenized equation which displays  explicit effective coefficients. The proof relies on the energetic method applied in the control-zone context. No critical thickness of the layers appears and the intersections of the layers have no distinct influence upon the asymptotic behaviour of the temperature.

The homogenization process in the  gridwork   case is worked out in Section~\ref{s:5}.  We obtain  the homogenized equation and explicit effective coefficients in this case also. The proofs are only
sketched because  the arguments follow the same lines as in  Section~\ref{s:4}. 

The significant difference  between  the effective coefficients obtained in the two  cases studied  here shows how important is the internal geometry of the vanishing superconductive material.

\section{Setting of the problem}\lab{s:2}

Let $\Om = I^3$, $I = \dis ]-\frac{1}{2},\frac{1}{2}[$ and $n\in\N$. From now on we use the notations
\be
\lab{1}
\veps = \frac{1}{2n+1},\quad \Z_\veps =\{k\in\Z,\quad \vert k\vert \leq n\}\quad\mbox{and}\quad 
I^k_\veps = \veps k+\veps I,\quad k\in\Z_\veps.
\ee 

Obviously, $\card\,\Z_\veps = {1}/{\veps}$ and $x\in\ov I$ if and only if there exists $k\in\Z_\veps$ such that $x\in\ov I^k_\veps$.

For any $i\in\{1,2,3\}$, we consider  ${r^i_\veps} \geq 0$, ${r^i_\veps} << \veps$, that is $\dis\frac{{r^i_\veps}}{\veps}\to 0$ when $\veps\to 0$. For any $k \in\Z_\veps$  we define: 
\be
\lab{2}
T^i_{\veps,k} = \{ x= (x_1,x_2,x_3)\in \Om,\;\vert x_i-\veps k\vert < {r^i_\veps}\},\;
{T^i_\veps} = \cup_{k\in\Z_\veps} T^i_{\veps,k}
\ee
\be
\lab{3}
T_\veps= \cup_{i=1}^3 T^i_\veps,\quad
T^{ij}_\veps = T^i_\veps\cap T^j_\veps,\quad\mbox{for}\quad 1\leq i < j\leq 3.
\ee

We consider that $\Om$ is occupied  by two materials  with highly  different conductivities; one is  highly conductive  and it is concentrated in the vanishing domain  of the thin layers $T_\veps$  ($\vert T_\veps\vert \to 0$) and the other forms the   ambiental  phase which has   conductivity of unity order. 

Assuming that we are given a source term $f_\veps$ in $\Om$ and some conductivity coefficients $a,b>0$, we consider the heat conduction problem, that is we study the temperature  $u_\veps$ which satisfies  in some sense  the system
\be
\lab{4}
-a\Delta  u_\veps = f_\veps\quad\mbox{in}\quad \Om\setminus T_\veps
\ee
\be
\lab{5}
-\frac{b}{\vert T_\veps\vert}\Delta  u_\veps = f_\veps\quad\mbox{in}\quad T_\veps
\ee
\be
\lab{6}
u_\veps=0\quad\mbox{on}\quad\pa\Om,
\ee
together with the natural transmission conditions on the interface $\pa T_\veps\setminus\pa\Om$.

More precisely,  assuming that $f_\veps\in H^{-1}(\Om)$, the variational  formulation of our  problem  is :

\vspace{0,3cm}
To find  $u_\veps\in H^1_0(\Om)$ such that 
\be
\lab{7}
a\int_{\Om\setminus T_\veps} \na u_\veps \na v \,+\, b\intb_{T_\veps} \na u_\veps\na v = \la f_\veps, v\ra,\quad \forall v\in H^1_0(\Om),
\ee
where we have used the notation
\be
\lab{8}
\intb_E\cdot \,= \frac{1}{\vert E\vert}\int_E\cdot\;\quad\mbox{for any measurable $E\subset\Om$}
\ee
and $\la\cdot,\cdot\ra$  is  the  duality product  between $H^{-1}(\Om)$ and $H^1_0(\Om)$.

Applying the Lax-Milgram theorem we obtain:

\bprop
\lab{p:2.1}
The variational equation (\ref{7}) has a unique solution $u_\veps\in H^1_0(\Om)$.
\eprop

The aim of this paper is to describe the asymptotic behaviour of the temperature $u_\veps$ as $\veps\to 0$,  assuming that the source term is weakly convergent:
\be
\lab{9}
\exists f\in H^{-1}(\Om)\quad\mbox{such that}\quad f_\veps\wto f\quad\mbox{in}\quad H^{-1}(\Om).
\ee

We have to remark here that   the boundedness of $(f_\veps)_\veps$ implies immediately:

\bprop
\lab{p:3.1}
$(u_\veps)_\veps$ is bounded in $H^1_0(\Om)$ and there exists $C>0$, independent of $\veps$, such that
\be
\lab{10}
\intb_{T_\veps}\vert\na u_\veps\vert^2\leq C.
\ee
\eprop

\section{The control-zone homogenization method}\lab{s:3}

In order to obtain further  results, specific  to the thin substructure considered here, we have to  introduce the following operators:

\bd
\lab{3.2}
To  any $u\in H^1(\Om)$   we associate  $G_\veps^i(u)\in L^2(\Om)$  defined by the following:
\be
\lab{11}
G_\veps^i(u)(x_1,x_2,x_3) = \sum_{k\in\Z_\veps}G^i_{\veps,k}(u)(\ov x_i)1_{I^k_\veps}(x_i),\quad
\ov x_i = (\cdots,\not{x}_i,\cdots)\in I^2
\ee
\be
\lab{12}
G^i_{\veps,k}(u) = \frac{1}{2} u\left\vert_{x_i=\veps k-{r^i_\veps}}\right. + 
\frac{1}{2} u\left\vert_{x_i=\veps k+{r^i_\veps}}\right.,\quad k\in\Z_\veps
\ee
where $\dis u\left\vert_{x_i=\veps k\pm{r^i_\veps}}\right. $ are the traces of $u$ on the corresponding faces of $T^i_{\veps,k}$.
\ed

These operators  have three basic properties, which were already proved in  \cite{GP08}:

\bprop
\lab{p:3.3}
For any $u\in H^1(\Om)$ and $i\in\{1,2,3\}$, we have 
\be
\lab{13}
\intb_{{T^i_\veps}} \vert G_\veps^i(u)-u\vert^2\leq{r^i_\veps}\left\vert\frac{\pa u}{\pa x_i}\right\vert_\Om^2
\ee
\be
\lab{14}
\intb_{{T^i_\veps} }\vert G_\veps^i(u)\vert^2= \vert G_\veps^i(u)\vert^2_\Om
\ee
\be
\lab{15}
 \vert G_\veps^i(u)-u\vert_\Om \leq\veps\left\vert\frac{\pa u}{\pa x_i}\right\vert_\Om
\ee
where $\vert\cdot\vert_\Om$ is denoting the norm of $L^2(\Om)$.
\eprop

The first  important  consequence is

\bthm
\lab{th:3.4}
There exists $C>0$, independent of $\veps$, such that 
\be
\lab{16}
\intb_{T_\veps}\vert u\vert^2\leq C\vert\na u\vert^2_\Om,\quad\forall u\in H^1_0(\Om).
\ee
\ethm
\bP
Let us notice that
$$
 \intb_{T_\veps}\vert u\vert^2\leq \sum_{i=1}^2\intb_{{T^i_\veps}}\vert u\vert^2\leq
2 \sum_{i=1}^2 \intb_{{T^i_\veps}}\left(\vert G_\veps^i(u)-u\vert^2+\vert G_\veps^i(u)\vert^2\right).
$$
Using (\ref{13})--(\ref{14}), it yields 
$$
\intb_{T_\veps}\vert u\vert^2\leq 2\,(\max_{i=1,2,3} r^i_\veps)\, \vert\na u\vert^2_\Om + 4\sum_{i=1}^2
\left(\vert G_\veps^i(u)-u\vert^2_\Om+\vert u\vert^2_\Om\right)
$$
and the proof is completed by (\ref{15}) and the Poincar\'e-Friedrichs inequality in $\Om$.
\eP

As  the techniques of  Section~3~\cite{BGPprep} can be used  to the domain $T^{ij}_\veps$ ($i < j$), then, according to \cite{BellieudB}, we have:

\bthm
\lab{th:3.5}
There exists $C>0$, independent of $\veps$, such that 
\be
\lab{17}
\intb_{T^{ij}_\veps} \vert u\vert^2\leq C\max\!\left(1,\,\veps^2\ln\frac{1}{r^i_\veps},\,\veps^2\ln\frac{1}{r^j_\veps} \right)\,\vert\na u\vert_\Om^2,\quad \forall u\in H^1_0(\Om).
\ee
\ethm

Finally, we remind an estimation of the same type proved in \cite{BGPprep} and which is associated to 
\be
\lab{17b}
T^0_\veps = \cap_{i=1}^3T^i_\veps.
\ee 

\bthm
\lab{th:3.6}
There exists $C>0$, independent of $\veps$, such that
\be
\lab{17c}
\intb_{T^0_\veps} \vert u\vert^2\leq C\max\!\left(1,\,\frac{\veps^3}{r_\veps} \right)\,\vert\na u\vert_\Om^2,\quad \forall u\in H^1_0(\Om).
\ee
\ethm

As the  vanishing highly conductive layers have thicknesses of size ${r^i_\veps}$, $0< {r^i_\veps} <<\veps$, we find the asymptotic behaviour of $u_\veps$  by applying a  control-zone method, that is  an energetic method using test-functions associated to an adequate  control-zone: a vanishing set which includes  the layers and has   much larger  thicknesses   ${R^i_\veps}$, where ${r^i_\veps} <<{R^i_\veps} << \veps$.

First, let us introduce the tools of this method. Denoting 
\be
\lab{25}
\RR= \{({R_\veps})_\veps,\; {r_\veps} << {R_\veps} << \veps\},
\ee
then, for any $({R^i_\veps})_\veps\in\RR$  we define the control-zone  of the present problem by
\be
\lab{26}
\C_\veps = \cup_{i=1}^3\C^i_\veps ,\quad 
\C^i_\veps= \cup_{k\in\Z_\veps}\C^i_{\veps,k},\quad i\in\{1,2,3\},
\ee
where for any $k\in\Z_\veps$ and $ i\in\{1,2,3\}$ we have 
\be
\lab{27}
\C^i_{\veps,k} = \{ x = (x_1,x_2,x_3)\in\Om,\quad \vert x_i-\veps k\vert < {R^i_\veps}\}.
\ee

The test-functions associated to this control-zone are defined  by  using the following    capacitary functions 
$w^i_\veps\in W^{1,\infty}(\Om)$ ($i =  1,2,3$), given by
\be
\lab{30}
w^i_{\veps}(x_1,x_2,x_3) = \left\{\ba{lcl} \dis 1-\frac{{r^i_\veps}}{{R^i_\veps}} & \mbox{if} & x = (x_1,x_2,x_3)\in {T^i_\veps} \\
1-\dis\frac{\vert x_i-\veps k\vert}{{R^i_\veps}} & \mbox{if} & x\in (\C^i_{\veps,k}\setminus T^i_{\veps,k})\;\mbox{for some $k\in\Z_\veps$} \\
0 & \mbox{if} & x\in\Om\setminus \C^i_\veps,\ea\right.
\ee
and the step  approximation operators  introduced  by

\bd
\lab{d:4.1}
To any $\varphi\in \DD(\Om)$ we associate $\varphi^i_\veps \in L^\infty(\Om)$  ($i=1,2,3$), as follows
\be
\lab{31}
\varphi^i_\veps(x_1,x_2,x_3) =\sum_{k\in\Z_\veps}  \varphi\left\vert_{x_i=\veps k}(\ov x_i) \right.1_{I^k_{{R^i_\veps}}}(x_i),
\ee
where  $I^k_{{R^i_\veps}} := \veps k + 2{R^i_\veps} I$. 
\ed

These operators have the following basic properties:
\be
\lab{32}
\vert\na w^i_\veps\vert_{\C_\veps}\leq \left(\frac{2}{\veps R^i_\veps}\right)^{1/2}
\ee

\be
\lab{33}
\vert\varphi-\varphi^i_\veps\vert_{L^\infty(\C^i_\veps)}\leq {R^i_\veps}\vert\na\varphi\vert_{L^\infty(\Om)}.
\ee

\be
\lab{34}
\vert\na\varphi^i_\veps\vert_{\C_\veps}\leq \left(\frac{2 R^i_\veps}{\veps}\right)^{1/2}
\vert\na\varphi\vert_{L^\infty(\Om)}.
\ee

\section{The reticulated case}\lab{s:4}

For any $i\in\{1,2,3\}$  let us denote $m_i\geq 0$ as the limit of 
\be
\lab{9b}
m_i = \lim_{\veps\to 0} \frac{\vert T^i_\veps\vert}{\vert T_\veps\vert}.
\ee
Obviously, we have
$$
m_1 + m_2 + m_3 = 1.
$$

In this section, we consider the case when 
\be
\lab{9c}
m_i >0,\quad\forall i\in\{1,2,3\},
\ee
that is the case when all the three parameters $r^i_\veps$ have the same order of magnitude with respect to $\veps$. This geometry  is called sometimes as the box-structure case.

We can present now the  preliminary convergence results:

\bprop
\lab{3.6}
There exists $u\in H^1_0(\Om)$ such that, on some subsequence, there hold
\be\lab{18}
u_\veps\wto u\quad\mbox{in}\quad H^1_0(\Om)
\ee
\be
\lab{19}
G_\veps^i(u_\veps)\to u\quad\mbox{in}\quad L^2(\Om),\quad\forall i\in \{1,2,3\}
\ee
\be
\lab{20}
\intb_{T_\veps} u_\veps v\to\int_\Om uv,\quad\forall v\in H^1_0(\Om).
\ee
Moreover, for any $i\in\{1,2,3\}$ we have:
\be
\lab{21}
\intb_{{T^i_\veps}}\frac{\pa u_\veps}{\pa x_j}v\to\int_\Om \frac{\pa u}{\pa x_j}v,\quad\forall v\in H^1_0(\Om),\quad \forall j\in\{1,2,3\},\quad j\not=i.
\ee
\eprop
\bP
The weak convergence (\ref{18}) follows from Proposition~\ref{p:3.1}. As $H^1_0(\Om)$ is  compactly embedded in $L^2(\Om)$, the strong convergence (\ref{19}) is obtained by using (\ref{15}) and (\ref{18}).

In order to prove (\ref{20}) we remark that 
\be
\lab{23}
\intb_{T_\veps} u_\veps v =  \frac{\vert T^0_\veps\vert}{\vert T_\veps\vert}\intb_{T^0_\veps} u_\veps v
-\sum_{1\leq i < j\leq 3} \frac{\vert T^{ij}_\veps\vert}{\vert T_\veps\vert}\intb_{T^{ij}_\veps} u_\veps v
+ \sum_{i=1}^3\frac{\vert {T^i_\veps}\vert}{\vert T_\veps\vert}\intb_{{T^i_\veps}}u_\veps v.
\ee
Using (\ref{10}), (\ref{17})  and (\ref{17c}), we prove that the first  two terms of (\ref{23}) are  converging to zero, as follows:
$$
\left\vert\frac{\vert T^0_\veps\vert}{\vert T_\veps\vert}\intb_{T^0_\veps} u_\veps v\right\vert
\leq C \frac{\vert T^0_\veps\vert}{\vert T_\veps\vert}
\max\!\left(1,\,\frac{\veps^{3/2}}{r^{1/2}_\veps} \right) \vert\na v\vert_\Om\to 0,
$$

$$
\left\vert\frac{\vert T^{ij}_\veps\vert}{\vert T_\veps\vert}\intb_{T^{ij}_\veps} u_\veps v\right\vert
\leq C \frac{\vert T^{ij}_\veps\vert}{\vert T_\veps\vert}
\max\!\left(1,\,\veps\ln^{\frac{1}{2}}\!\!\frac{1}{r^i_\veps},\,\veps\ln^{\frac{1}{2}}\!\!\frac{1}{r^j_\veps} \right) \vert\na v\vert_\Om\to 0.
$$
Taking into account  (\ref{9b}), the proof of (\ref{20})  is completed by applying to  Proposition~2.15 of  \cite{GP08}, where we have proved that
\be
\lab{24}
\intb_{T^i_\veps} u_\veps v\to\int_\Om uv,\quad\forall v\in H^1_0(\Om).
\ee

Next, let $\varphi\in\DD(\Om)$ and $i\in\{1,2,3\}$. Denoting  by $\nu = (\nu_1,\nu_2,\nu_3)$ the outward normal to $\pa T_\veps$ we obviously have $\varphi\nu_j = 0$ on $\pa {T^i_\veps}$ ($j\not=i$) and hence:
$$
 \intb_{{T^i_\veps}}\frac{\pa u_\veps}{\pa x_j}\varphi= - \intb_{{T^i_\veps}}u_\veps\frac{\pa\varphi}{\pa x_j}.
$$
Using (\ref{24}), from the previous relation it follows:
$$
 \intb_{{T^i_\veps}}\frac{\pa u_\veps}{\pa x_j}\varphi \to -\int_\Om u\frac{\pa\varphi}{\pa x_j} = \int_\Om\frac{\pa u}{\pa x_j}\varphi.
$$
The proof of (\ref{21}) is completed by continuity, using (\ref{10}) and (\ref{16}).
\eP

Now  we can present our main result. 

\bthm
\lab{th:4.2}
$(u_\veps)_\veps$ is weakly convergent in $H^1_0(\Om)$. Its limit, $u\in H^1_0(\Om)$, is the only solution of the equation
\be
\lab{35}
-\sum_{i=1}^3\left(a+\frac{b}{3}(1-m_i)\right)\frac{\pa^2u}{\pa x_i^2} = f\quad\mbox{in}\quad\Om.
\ee
\ethm
\bP
For any $i\in\{1,2,3\}$, let $(\R^i_\veps)_\veps\in\RR$ and $\varphi\in\DD(\Om)$. Using the definitions (\ref{30})--(\ref{31}), we denote 
\be
\lab{36}
v_\veps(\varphi) =\sum_{i=1}^3\left(\left(1-\frac{r^i_\veps}{R^i_\veps}\right)\varphi + (\varphi^i_\veps-\varphi) w^i_\veps\right).
\ee
We notice that 
\be
\lab{37}
\!\!\!
\varphi^i_\veps(x)w^i_\veps(x) = \left\{\ba{l}
\dis\left(1-\frac{r^i_\veps}{R^i_\veps}\right)\varphi\vert_{x_i=\veps k}(\ov x_i) ,\; x\in T^i_{\veps,k}\quad\mbox{for some}\quad k\in\Z_\veps\\
\\
\dis\left(1-\frac{\vert x_i-\veps k\vert}{R^i_\veps}\right)\varphi\vert_{x_i=\veps k}(\ov x_i) , \; x\in (\C^i_{\veps,k}\setminus T^i_{\veps,k}),\; k\in\Z_\veps\\
\\
0,\; x\in\Om\setminus\C^i_\veps\ea\right.
\ee
and hence $v_\veps(\varphi)\in W^{1,\infty}_0(\Om)$. Then we set $v = v_\veps(\varphi)$ in (\ref{7}) and it follows 
\be
\lab{38}
\ba{l}
a\dis\sum_{i=1}^3\dis\left(1-\frac{r^i_\veps}{R^i_\veps}\right) \int_{\Om\setminus\C_\veps}\na u_\veps\na\varphi + a \int_{\C_\veps\setminus T_\veps} \na u_\veps\na v_\veps(\varphi)\, + \\
\quad\quad\quad\quad\quad
+\, b\dis\sum_{i=1}^3\dis\left(1-\frac{r^i_\veps}{R^i_\veps}\right) \intb_{T_\veps}\na u_\veps\na\varphi^i_\veps = 
\la f_\veps,v_\veps(\varphi)\ra.\ea
\ee
As $\vert\C_\veps\vert\to 0$ it follows that
\be
\lab{39}
\left(1-\frac{r^i_\veps}{R^i_\veps}\right)(\na\varphi)1_{\Om\setminus\C_\veps}\to \na\varphi\quad\mbox{strongly in}\quad L^2(\Om)
\ee
and using (\ref{18}) we obtain the convergence of the first left-hand side term of (\ref{38}):
\be
\lab{40}
a\sum_{i=1}^3\dis\left(1-\frac{r^i_\veps}{R^i_\veps}\right) \int_{\Om\setminus\C_\veps}\na u_\veps\na\varphi
\, \to \, 3 a\! \int_\Om\na u\na\varphi.
\ee

The second left-hand side term of (\ref{38}) is converging to zero. Indeed, as $(u_\veps)_\veps$ is bounded in $H^1_0(\Om)$, we have 
\be
\lab{41}
\left\vert \;a \int_{\C_\veps\setminus T_\veps} \na u_\veps\na v_\veps(\varphi)\,\right\vert\leq C
\vert\na v_\veps(\varphi)\vert_{\C_\veps}.
\ee
Moreover, by a straightforward computation we find that 
\be
\lab{42}
\vert\na v_\veps(\varphi)\vert_{\C_\veps}\leq C\vert\na\varphi\vert_{\C_\veps} + C
\sum_{i=1}^3 \left(
\vert\varphi^i_\veps - \varphi\vert_{L^\infty(\C^i_\veps)} \vert\na w^i_\veps\vert_{\C_\veps} +
\vert\na_{\ov x_i} \varphi^i_\veps\vert_{\C_\veps} \right)
\ee
and the assertion is proved by using (\ref{32})--(\ref{34}):
\be
\lab{43}
\vert\na v_\veps(\varphi)\vert_{\C_\veps}\leq C\vert\na\varphi\vert_{L^\infty(\Om)} \sum_{i=1}^3\left(\frac{R^i_\veps}{\veps}\right)^{1/2}  \to \,\, 0.
\ee

For the convergence of the third left-hand side term of (\ref{38}) we notice that 
\be
\lab{44}
\ba{l}
\!\!\!
\dis\intb_{T_\veps}\na u_\veps\na\varphi^i_\veps  =\dis \frac{1}{\vert T_\veps\vert}
\int_{T^i_\veps}\na u_\veps\na\varphi^i_\veps = \\
=\dis
\frac{1}{\vert T_\veps\vert}\sum_{k\in\Z_\veps} \int_{T^i_{\veps,k}}\na_{\ov x_i}u_\veps(x)
\na_{\ov x_i}\varphi\vert_{x_i=\veps k}(\ov x_i)\;dx\ea.
\ee
Using (\ref{10}) and the smoothness of $\varphi$ we have
$$
\left\vert\frac{1}{\vert T_\veps\vert}\sum_{k\in\Z_\veps} \int_{T^i_{\veps,k}}\na_{\ov x_i}u_\veps(x)
\left(\na_{\ov x_i}\varphi\vert_{x_i=\veps k}(\ov x_i) - \na_{\ov x_i}\varphi(x) \right)\;dx\right\vert\leq
$$
$$
\leq C(\varphi)r_\veps\frac{1}{\vert T_\veps\vert}\int_{T^i_\veps}\vert\na_{\ov x_i}u_\veps\vert\leq C(\varphi) r_\veps \,\to \, 0
$$
and thus  we replace (\ref{44}) by
$$
\lim_{\veps\to 0}\intb_{T_\veps}\na u_\veps\na\varphi^i_\veps = \lim_{\veps\to 0}
\frac{\vert T^i_\veps\vert}{\vert T_\veps\vert} \int_{T^i_{\veps}}\na_{\ov x_i}u_\veps
\na_{\ov x_i}\varphi.
$$
Next, recalling (\ref{9b}) and applying (\ref{21}) we obtain
\be
\lab{45}
\ba{l}
b\dis \sum_{i=1}^3\dis\left(1-\frac{r^i_\veps}{R^i_\veps}\right) \intb_{T_\veps}\na u_\veps\na\varphi^i_\veps\to
\dis b\sum_{i=1}^3m_i \int_\Om\na_{\ov x_i} u\na_{\ov x_i}\varphi = \\
\quad\quad\quad\quad\quad\quad\quad\quad
\quad\quad\quad\quad\quad\quad
=b\dis \sum_{i=1}^3(1-m_i)\int_\Om\frac{\pa u}{\pa x_i}\frac{\pa\varphi}{\pa x_i}.\ea
\ee

Finally, noticing that (\ref{32}) and (\ref{33}) imply 
$$
\left\vert\la f_\veps,(\varphi^i_\veps-\varphi)w^i_\veps\ra\right\vert\leq C\vert w^i_\veps\vert_{\C_\veps}
\vert\na\varphi\vert_{L^\infty(\Om)} +
C\vert\na w^i_\veps\vert_{\C_\veps}\vert\varphi^i_\veps-\varphi\vert_{L^\infty(\C^i_\veps)}\leq
$$
$$
\leq C(\varphi)\left( \vert\C_\veps\vert^{1/2} + \left\vert\frac{R^i_\veps}{\veps}\right\vert^{1/2}\right)\,\to\, 0,
$$
we obtain the convergence of the right-hand side term of (\ref{38}), that is:
\be
\lab{46}
\la f_\veps, v_\veps(\varphi)\ra = \sum_{i=1}^3\dis\left(1-\frac{r^i_\veps}{R^i_\veps}\right)\la f_\veps, \varphi\ra +
\sum_{i=1}^3\la f_\veps,(\varphi^i_\veps-\varphi)w^i_\veps\ra \,\to\, 3\,\la f,\varphi\ra.
\ee
Resuming, we can say that if we pass (\ref{38})  to the limit, then we obtain:
\be
\lab{47}
3a\int_\Om\na u\na\varphi + b\sum_{i=1}^3(1-m_i)\int_\Om\frac{\pa u}{\pa x_i}\frac{\pa\varphi}{\pa x_i} =
3\,\la f,\varphi\ra,\quad\forall\varphi\in\DD(\Om),
\ee
which is, by continuity, the variational formulation of (\ref{35}). The convergence of the entire sequence $(u_\veps)_\veps$ follows from the unicity of the solution of (\ref{35}).
\eP

\brem
\lab{4.3}
In spite of the vanishing volume of the rectangular honeycomb structure  the homogenized behaviour is anisotropic, except the case when $m_1=m_2=m_3 = \dis\frac{1}{3}$. 
\erem

\brem
\lab{4.4}
The reticulated domain occupied by all the intersections of the layers has no distinct influence upon the asymptotic distribution of the temperature. For $i<j$ and $k\in\{i,j\}$, this follows from 
\be
\lab{48}
\!\!\!
\ba{l}
\dis\left\vert\frac{1}{\vert T_\veps\vert}\int_{T^{ij}_\veps}\na u_\veps\na\varphi^k_\veps\right\vert\leq C\vert\na\varphi\vert_{L^\infty(\Om)}\left(\frac{1}{\vert T_\veps\vert}\int_{T^{ij}_\veps}\vert\na u_\veps\vert^2\right)^{\frac{1}{2}}
\left(\frac{\vert T^{ij}_\veps\vert}{\vert T_\veps\vert}\right)^{\frac{1}{2}}\leq \\
\quad\quad\quad\quad\quad\quad\quad\quad\quad\quad\quad\quad\quad\quad
\dis\leq C(\varphi)\left(\frac{r_\veps}{\veps}\right)^{\frac{1}{2}} \,\to\, 0,\ea
\ee
where we have used (\ref{10}).
\erem

\section{The gridwork case}\lab{s:5}

Using the notations and definitions of the previous sections, the gridwork case corresponds to the situation when the horizontal layers are missing and the other two layers have the same order of magnitude with respect to $\veps$. It follows that 
\be
\lab{51}
T_\veps = T^1_\veps\cup T^2_\veps
\ee
and there exist 
\be
\lab{52}
m_i = \lim_{\veps\to 0}\frac{\vert T^i_\veps\vert}{\vert T_\veps\vert}\quad (i=1,2),\quad\mbox{such that}\quad m_1+m_2 = 1.
\ee
For consistency with the honeycomb case, we also define here 
\be
\lab{53}
m_3 = 0.
\ee

Obviously, for $i\in\{1,2\}$ th properties of the corresponding operators defined in  Section~\ref{s:3} are still valid. Then we prove results similar to those of Proposition~\ref{3.6}, with only one mentionable distinction:

For any $i\in\{1,2\}$ we have (\ref{21}). 

Consequently, the homogenization result in this case is the following:

\bthm
\lab{5.1}
$(u_\veps)_\veps$ is weakly convergent in $H^1_0(\Om)$. Its limit, $u\in H^1_0(\Om)$, is the only solution of the equation:
\be
\lab{54}
-\sum_{i=1}^3\left(a+\frac{b}{2}(1-m_i)\right)\frac{\pa^2u}{\pa x_i^2} = f\quad\mbox{in}\quad\Om.
\ee
\ethm
\bP
The proof is completely analogous to that of Theorem~\ref{th:4.2}, with the difference that the control-zones are oriented in only two directions. That is, with the test-function
\be
\lab{55}
v_\veps(\varphi) =\sum_{i=1}^2\left(\left(1-\frac{r^i_\veps}{R^i_\veps}\right)\varphi + (\varphi^i_\veps-\varphi) w^i_\veps\right)
\ee
we get the same convergence as in (\ref{45}): 
\be
\lab{56}
\ba{l}
b\dis \sum_{i=1}^2\dis\left(1-\frac{r^i_\veps}{R^i_\veps}\right) \intb_{T_\veps}\na u_\veps\na\varphi^i_\veps\to
\dis b\sum_{i=1}^2m_i \int_\Om\na_{\ov x_i} u\na_{\ov x_i}\varphi = \\
\quad\quad\quad\quad\quad\quad\quad\quad
\quad\quad\quad\quad\quad\quad
=b\dis \sum_{i=1}^3(1-m_i)\int_\Om\frac{\pa u}{\pa x_i}\frac{\pa\varphi}{\pa x_i},\ea
\ee
where we have taken in account the convention (\ref{53}) also. From here the proof  can be completed identically. 
\eP

\brem
\lab{5.2}
Assuming that the quantity $\vert T_\veps\vert$ is the same in the both cases that we considered, the significant difference between the equations (\ref{35}) and (\ref{55}) shows how important is the internal geometry of the vanishing superconductive material. 
\erem

\vs5

{\bf Acknowledgements.} This work has been accomplished during the visit of
 D.~Poli\c{s}evschi at the I.R.M.A.R.'s Department of Mechanics
(University of Rennes~1), whose support is gratefully acknowledged.


\begin{thebibliography}{1}

\bibitem{BakhvalovP89}
{\sc N.S.~Bakhvalov, G.P.~Panasenko},
\newblock Homogenization: averaging processes in periodic media, Math. Appl.  (Sov. Series), {\bf 36},
\newblock Kluwer Acad. Pub. Group, Dordrecht, 1989.


\bibitem{BellieudB}
{\sc M.~Bellieud, G.~Bouchitt\'e},
\newblock Homogenization of elliptic problems in a fiber reinforced structure.
  Non local effects,
\newblock {\em Ann. {S}cuola {N}orm. {S}up. {P}is {C}l. {S}ci.(4)}, {\bf 26(3)},
  1998, 407--436.


\bibitem{BGP06}
{\sc F.~Bentalha, I.~Gruais, D.~Polisevski},
\newblock Homogenization of a conductive suspension in a Stokes-Boussinesq flow,
\newblock {\em Applicable Analysis}, {\bf 85(6-7)}, 2006, 811--830.



\bibitem{BGP}
{\sc F.~Bentalha, I.~Gruais, D.~Polisevski},
\newblock Diffusion process in a  rarefied binary
structure,
\newblock {\em Revue Roumaine de Math\'ematiques Pures et Appliqu\'ees}, {\bf 52(1)}, 2007, 9--34.



\bibitem{BGPprep}
{\sc F.~Bentalha, I.~Gruais, D.~Polisevski},
\newblock {Diffusion in a highly rarefied binary structure of general periodic shape},
\newblock {\em Applicable Analysis}, {\bf 87(6)}, 2008, 635--655.


\bibitem{BouchitteP96}
{\sc G.~Bouchitt\'e, C.~Picard},
\newblock {Singular perturbations and homogenization in stratified media},
\newblock {\em Applicable Analysis}, {\bf 61(3-4)}, 1996,  307--341.

\bibitem{GP07}
{\sc I.~Gruais, D.~Polisevski},
\newblock{Homogenizing a critical binary structure of finite diffusivities}, 
\newblock{\em Asymptotic Analysis}, {\bf 55(1)}, 2007, 85--102.


\bibitem{GP08}
{\sc I.~Gruais, D.~Polisevski},
\newblock {Periodic structures separated by highly conductive thin layers},  
\newblock{\em Pr\'epublication~IRMAR}, {\bf 08-28}, 2008, 1--21.


\bibitem{HeronMP94}
{\sc B.~H{\'e}ron, J.~Mossino, C.~Picard},
\newblock {Homogenization of some quasilinear problems for stratified media with low and high conductivities},
\newblock {\em Differential Integral Equations}, {\bf 7(1)}, 1994, 157--178.



\bibitem{Murattartar85}
{\sc F.~Murat, L.~Tartar},
\newblock  {Calcul des variations et homog\'en\'eisation},
\newblock in {'' Homogenization methods: theory and applications in physics''}, {\em Collect. Dir. \'Etudes Rech. \'Elec. France}, {\bf 57}, Eyrolles, Paris, 1985, 319--370.




\bibitem{SJP99}
{\sc D.~Cioranescu, J.~Saint-Jean-Paulin},
\newblock  Homogenization of reticulated structures,
\newblock  Springer-Verlag, New-York, 1999.


\bibitem{ZhikovP06}
{\sc V.V.~Zhikov, S.E.~Pastukhova},
\newblock Derivation of the limit equations of elasticity theory of thin sets,
\newblock{\em Journal of Mathematical Sciences}, {\bf 135(1)}, 2006, 2637--2665.

 



\end{thebibliography}

%
\newcommand{\noopsort}[1]{}

\vs5

* Universit\'e de Rennes1, I.R.M.A.R, Campus de Beaulieu,
35042 Rennes Cedex (France)

\vs5

** I.M.A.R., P.O. Box 1-764, RO-014700  Bucharest (Romania).

\end{document}